\documentstyle[12pt]{article}
\newtheorem{theorem}{Theorem}[section]
\newtheorem{corollary}[theorem]{Corollary}
\newtheorem{definition}{Definition}

\newtheorem{remark}[theorem]{Remark}
\newtheorem{lemma}[theorem]{Lemma}
\newtheorem{proposition}[theorem]{Proposition}
\title{Orthocompactness and semi-stratifiability in the density
topology\thanks{1991 Math.\ Subject Classification --- Primary:
28A05, 54G20, 54-02; Secondary: 03E15, 04A15, 54A10.\newline Key
words and phrases --- density topology, measure, pseudocompact,
orthocompact, semi-stratifiable. \newline Research supported
partially by the Ella and Georg Ehrnrooth Foundation at Merita
Bank, Finland.}}
\author{Julian Dontchev\\Department of Mathematics\\University
of Helsinki\\PL 4, Yliopistonkatu 5\\00014 Helsinki\\Finland}
\date{}
\begin{document}
\baselineskip=20pt plus 1pt minus 1pt
\maketitle
\begin{abstract} 
The density topology $\cal T$ is a topology on the real line,
finer than the usual topology, having as its open sets the
measurable subsets of ${\bf R}$, which are of density 1 at
each of their points. The aim of this paper is to determine which
subsets of the density topology are semi-stratifiable,
orthocompact and weakly hereditarily pseudocompact.
\end{abstract}

\section{The density topology: introduction}\label{s1}

The density topology on the real line is a finer topology than
the usual Euclidean one, serving often as a counterexample
in modern Topology. The open sets of the density topology are the
measurable subsets of the real line, which are of density 1 at
each of their points. The density topology is closely related to
the measure structure of the real line and is a {\em bridge}
between Analysis and Topology. In fact, the density topology is
the place where Analysis, Measure theory, Set theory and Topology
celebrate their consequent happy reunion, Topology probably being
the host.

It was not until 1952, when the density topology was introduced
by Haupt and Pauc \cite{HP0}. For the first time the properties
of the density topology were studied in 1961 by Goffman and
Waterman \cite{GW1}. Two very important contributions on the
study are due to Tall \cite{T1,T2}. In connection with
approximate continuity Goffman together with Neugebauer and
Nishiura continued the investigation on the density topology.

The effort to explain the striking parallels between the theorems
concerning nullsets and first category sets has been for years
a major block of research in the field of real analysis. Indeed
the density topology on the real line is a topology where these
two concepts coincide. As the Lebesgue density plays a central
r\^{o}le in the study of real functions it has its significant
impact on the density topology.

The density topology, with its rich structure, is very often a
useful counterexample in the study of general topological spaces.
The aim of this paper is to investigate the density topology in
the light of the some popular (and less popular) concepts in
Modern Topology, namely, we try to determine which subsets of the
density topology are semi-stratifiable, orthocompact and weakly
hereditarily pseudocompact.

\begin{definition}\label{d1}
{\em A measurable set $E \subseteq {\bf R}$ has density $d$
at $x \in {\bf R}$ if $$\lim_{h \rightarrow 0} \frac{m(E \cap
[x-h,x+h])}{2h}$$ exists and is equal to $d$. Set $\phi(E) = \{
x \in {\bf R} \colon d(x,E) = 1 \}$. The open sets of the
density topology $\cal T$ are those measurable sets $E$ that
satisfy $E \subseteq \phi(E)$. Clearly, the density topology
$\cal T$ is finer than the usual topology on the real line.}
\end{definition}

Tall \cite{T1} proved that every subset of the density topology
is the union of a CCC-set and a closed discrete set. He also
proved that the density topology is hereditarily subparacompact
and hence hereditarily submetacompact (= hereditarily
$\theta$-refinable \cite{WoW1}) and hereditarily countably
subparacompact. In the same paper, Tall showed that the
${\aleph}_{1}$-compact subsets of the density topology are
hereditarily Lindel\"{o}f, and that all collectionwise Hausdorff
and all $\sigma$-metacompact subspaces of the density topology
are the union of a hereditarily Lindel\"{o}f and a closed
discrete set. Moreover, every countably paracompact subspace $A$
of the density topology with $|RO(A)| < {\kappa}$ is
$\kappa$-compact \cite{T1} and if $A \subseteq ({\bf R},{\cal
T})$ is countably paracompact, then $A$ is the union of a
$2^{{\aleph}_{0}}$-compact set and a closed discrete set.

Since the density topology is not even $\sigma$-metacompact (but
subparacompact) \cite{T1}, it is natural to ask if the density
topology is orthocompact or semi-stratifiable. We also determine
the subsets of the density topology which have some properties
in the vicinity of pseudocompactness.

\section{Semi-stratifiability in the density topology}\label{s2}

It is well-known that every semi-stratifiable space is
subparacompact. Since the density topology is hereditarily
subparacompact and since semi-stratifiability is a hereditary
property, it is natural to ask if the density topology is
(hereditarily) semi-stratifiable.

\begin{theorem}\label{semstr1}
For a subset $A$ of the density topology the following conditions
are equivalent:

{\rm (1)} $A$ is a nullset.

{\rm (2)} $A$ is semi-stratifiable.
\end{theorem}

{\em Proof.} (1) $\Rightarrow$ (2) Clearly every nullset of the
density topology is semi-stratifiable, since the nullsets are
precisely the closed and discrete subsets.

(2) $\Rightarrow$ (1) Assume that $A$ is semi-stratifiable and
that $A$ is not a nullset. Let $B$ be a subset of $A$ with the
smallest possible cardinality $\kappa$ such that $B$ is not a
nullset and accordingly well-order $B$. Note that due to
heredity, $B$ is also semi-stratifiable. For each $x \in B$ let
$W\{x\} = B \setminus B_x$, where $B_x = \{ y \colon y < x \}$.
Note that each $B_x$ is a nullset and hence closed in the density
topology. Moreover, if $x \in W\{y\}$, then $y \not\in W\{x\}$.
Thus $W$ is an antisymmetric neighbornet of $B$ (see \cite{J0}
for the definition of a neighbornet). By \cite[Corollary
4.9]{J0}, $B$ is $\sigma$-discrete. Since in the density topology
every discrete set is a nullset, then $B$ is a nullset. By
contradiction, $A$ is a nullset. $\Box$

\begin{corollary}
The density topology is not semi-stratifiable and the
stratifiable subsets are precisely the closed and discrete ones.
\end{corollary}

\section{$N$-pseudocompactness in the density topology}\label{s3}

It is well-known that closed subspaces of pseudocompact spaces
are not necessarily pseudocompact. If every closed (resp.\
nowhere dense) subspace of a pseudocompact space $(X,\tau)$ is
pseudocompact, then $X$ is called {\em weakly hereditarily
pseudocompact} (resp.\ {\em $N$-pseudocompact}). Clearly, every
countably compact space is weakly hereditarily pseudocompact and
within the class of nodec spaces (= nowhere dense subsets are
closed) weakly hereditarily pseudocompactness implies
$N$-pseudocompactness. The following result relates to Theorem
3.6 from \cite{T1}.

\begin{theorem}\label{tpse1}
For a subset $A$ of the density topology the following conditions
are equivalent:

{\rm (1)} $A$ is (weakly) hereditarily pseudocompact.

{\rm (2)} $A$ is $N$-pseudocompact.

{\rm (3)} $A$ is finite.
\end{theorem}

{\em Proof.} (1) $\Rightarrow$ (2) is obvious, since the density
topology is nodec.

(2) $\Rightarrow$ (3) In the notion of \cite[Theorem
3.6]{T1}, we need to prove that $A$ is countably compact. Since
$A$ is Tychonoff and pseudocompact, then $A$ is lightly compact,
i.e.\ every countable open cover has a finite dense subsystem.
Let ${\cal U} = \{ U_i \colon i \in I \}$ be a countable open
cover of the subspace $A$ and let $\{ U_i \colon i \in F \}$ be
a finite dense subsystem of $\cal U$ in $A$. If $W = A \setminus
\cup_{i \in F} U_i$ is nonempty, then $W$ is nowhere dense in $A$
and hence closed and discrete \cite{T1}. Since $A$ is
$N$-pseudocompact, then $W$ is pseudocompact and hence lightly
compact. Thus the discrete subspace $W$ must be finite. Hence,
$A$ is countably compact and so finite.

(3) $\Rightarrow$ (1) is obvious. $\Box$

\begin{remark}
{\em (i) By Theorem~\ref{tpse1} all weakly hereditarily lightly
compact subspaces of the density topology are finite. 

(ii) In Theorem~\ref{tpse1}, `$N$-pseudocompact' can not be
replace by `$Z$-pseudocompact'.}
\end{remark}

Topological spaces in which every locally finite collection of
open sets is countable are called {\em pseudo
$\aleph_{1}$-compact}.

\begin{proposition}\label{t05}
{\rm (i)} The density topology is pseudo $\aleph_{1}$-compact;

{\rm (ii)} The density topology is not $\aleph_{1}$-compact.
\end{proposition}

{\em Proof.} In order to observe (i), note that even point-finite
collections of open sets are countable provided the space is both
Baire and CCC as such is the case with the density topology
\cite{T1}. Since submetacompact $\aleph_{1}$-compact spaces must
be Lindel\"{o}f \cite{WoW1}, then (ii) is also clear. $\Box$

\begin{remark}
{\em (i) The density topology is even {\em strongly pseudo
$\aleph_{1}$-compact}, i.e.\ every point-finite collection of
open sets is countable;

(ii) The density topology is not a stable quasi-pseudo metric
space, since as proved in \cite{J6} in stable quasi-pseudo metric
spaces pseudo $\aleph_{1}$-compactness implies separability;

(iii) Recall that a topological space is called {\em mildly
countably compact} \cite{Staum1} if every disjoint open cover of
$X$ has only a finite number of non-empty members or equivalently
if there is no continuous function from $X$ onto the integers.
Clearly, the density topology is mildly countably compact.}
\end{remark}

\section{In the vicinity of orthocompactness}\label{s4}

A topological property strictly weaker than metacompactness is
orthocompactness. Recall first that a family $\cal U$ of open
subsets is called {\em interior preserving} if for every $\cal
V \subseteq \cal U$, $\cap \cal V$ is open. A topological space
$(X,\tau)$ is called {\em (countably) orthocompact} \cite{FL1}
if every (countable) open cover of $X$ has an interior-preserving
open refinement. Orthocompactness is a strictly weaker property
than metacompactness, since all linearly ordered topological
spaces are orthocompact. Also, all principal spaces (= union of
closed sets is always a closed set) are orthocompact and all
P-spaces are countably orthocompact.

Orthocompact spaces were first considered in 1966 by Sion and
Willmott \cite{SW1} as the spaces having the property Q. The name
orthocompact was given by Arens.

\begin{proposition}
The density topology is (hereditarily) countably orthocompact.
\end{proposition}

{\em Proof.}  As mentioned in \cite{D1}, every countably
metacompact (= countably submetacompact \cite{D1}) space is
countably orthocompact. Note that the density topology is even
countably subparacompact. $\Box$

Recall that a space $(X,\tau)$ is called {\em
$\sigma$-orthocompact} \cite{FL0} if $\cal C$ is an open cover
of $X$, then there exists an open refinement ${\cal R} =
\cup_{n=1}^{\infty} {\cal R}_{n}$ of $\cal C$ such that for each
$n \in \omega$, ${\cal R}_{n}$ is interior preserving.

\begin{theorem}\label{tdo1}
{\rm (CH)} The density topology has a dense hereditarily
orthocompact subspace of power continuum.
\end{theorem}

{\em Proof.} First, we establish that the density topology has
a dense $\sigma$-orthocompact subspace of power continuum. It is
proved in \cite{DTW1} that every CCC, Baire, dense-in-itself
space with $\pi$-weight $\leq 2^{{\aleph}_{0}}$ contains a dense
generalized Lusin subspace (= every nowhere dense subset has
cardinality $\leq 2^{{\aleph}_{0}}$) of power continuum. Since
the density topology is cometrizable \cite{T2}, then its
$\pi$-weight $\leq 2^{{\aleph}_{0}}$. Let $S$ be the dense
subspace (of the density topology) in question. By \cite[Theorem
3.1]{T1}, $S = A \cup B$, where $A$ is CCC and $B$ is nowhere
dense. Clearly, $B$ is countable. Thus in order to show that $S$
is $\sigma$-orthocompact, it suffices to verify that $A$ is
$\sigma$-orthocompact. Let $\cal U$ be an open cover of $A$. Let
${\cal V} = \{ V_n \colon n \in \omega \}$ be a maximal disjoint
collection of open sets such that each one is included in a
member of $\cal U$. Since $A \setminus \cup {\cal V}$ is nowhere
dense in $A$, then it is nowhere dense in $S$ and hence
countable. Since $\cal V$ is interior preserving, then $A$ is
$\sigma$-orthocompact. Thus $S$ is $\sigma$-orthocompact and
since $S$ is also perfect due to heredity, then $S$ is
hereditarily orthocompact. $\Box$

Recall that a topological space $(X,\tau)$ is called {\em weakly
orthocompact} (in the sense of Peregudov) \cite{Per1} if every
open cover $\cal U$ of $X$ has an open refinement $\cal V$ such
that for each $x \in X$ the set $\cap {\cal V}_{x}$ has nonempty
interior. In \cite{Scott1}, Scott defined another form of
orthocompactness, which he also called {\em weakly orthocompact}.
Scott's definition requires that directed open covers have
interior preserving open refinements. A topological space
$(X,\tau)$ is called {\em semi-metacompact} \cite{FL0,Per1} if
every open cover of $X$ has an open-finite refinement, where an
open-finite cover means that no nonempty open set is a subset of
infinitely many members of the cover \cite{Per1}.

\begin{proposition}
{\rm (CH) (i)} The density topology is a weakly orthocompact
space in the sense of Peregudov.

{\rm (CH) (ii)} The density topology is not semi-metacompact.
\end{proposition}

{\em Proof.} (i) Follows from \cite[Theorem 3]{Per1} and the fact
that under CH, the density topology is meta-Lindel\"{o}f.

(ii) Follows from the fact that all semi-metacompact weakly
orthocompact spaces are metacompact. $\Box$

A slightly stronger form of orthocompactness was recently
considered by Junnila and K\"{u}nzi. A topological space
$(X,\tau)$ is called {\em ortho-refinable} \cite{J5} provided
that for each open cover $\cal C$ of $X$ there are an ordinal
$\delta$ and a decreasing chain $(T_{\alpha})_{{\alpha} <
{\delta}}$ of transitive partial neighbornets on $X$ so that for
each $x \in X$ there exists $\alpha < \delta$ such that $x \in
{\rm St} (x,{\cal T}_{\alpha}) \subseteq C$ for some $C \in {\cal
C}$. Here ${\cal T}_{\alpha} = \{ T_{\alpha} (x) \colon x \in
T_{\alpha} (X) \}$. Concerning neighbornets the reader may refer
to \cite{J0}. The class of ortho-refinable spaces is placed
between the classes of spaces having ortho-bases and
orthocompactness \cite{J5}.

\begin{lemma}\label{JK2}
{\em \cite{J5}} A submetacompact space is ortho-refinable if and
only if it is orthocompact.
\end{lemma}

\begin{theorem}\label{tdo2}
{\rm (CH)} The density topology has a dense hereditarily
ortho-refinable subspace of power continuum.
\end{theorem}

{\em Proof.} Follows from Theorem~\ref{tdo1} and Lemma~\ref{JK2}.
$\Box$

In the notion of Lemma~\ref{JK2}, the orthocompactness of the
density topology would imply its ortho-refinability. On the other
hand, if the density topology is a $\sigma$-orthocompact space,
then it clearly would be hereditarily orthocompact. 

{\bf Question.} Is the density topology ($\sigma$-)orthocompact?

\
E-mail: {\tt dontchev@cc.helsinki.fi}, {\tt
dontchev@e-math.ams.org}
\
\end{document}